\documentclass[11pt]{article}
\usepackage{amsmath,amssymb} 

\begin{document}
\title{Exceptional Dehn surgery on large arborescent knots}
\author{Ying-Qing Wu}
\date{}
\maketitle

\footnotetext[1]{ Mathematics subject classification (1991):  {\em Primary 
57N10.}}
\footnotetext[2]{ Keywords and phrases: Arborescent knots, exceptional Dehn surgery}

\begin{abstract}
  A Dehn surgery on a knot $K$ in $S^3$ is exceptional if it produces
  a reducible, toroidal or Seifert fibred manifold.  It is known that
  a large arborescent knot admits no such surgery unless it is a type
  II arborescent knot.  The main theorem of this paper shows that up
  to isotopy there are exactly three large arborescent knots admitting
  exceptional surgery, each of which admits exactly one exceptional
  surgery, producing a toroidal manifold.
\end{abstract}

\newcommand{\proof}{\noindent {\bf Proof.} }
\newcommand{\qed}{\quad $\Box$ \medskip}

\newtheorem{thm}{Theorem}[section]
\newtheorem{prop}[thm]{Proposition} 
\newtheorem{lemma}[thm]{Lemma} 
\newtheorem{cor}[thm]{Corollary} 
\newtheorem{defn}[thm]{Definition} 
\newtheorem{notation}[thm]{Notation} 
\newtheorem{qtn}[thm]{Question} 
\newtheorem{example}[thm]{Example} 
\newtheorem{remark}[thm]{Remark} 
\newtheorem{conj}[thm]{Conjecture} 
\newtheorem{prob}[thm]{Problem} 
\newtheorem{rem}[thm]{Remark} 

\newcommand{\bdd}{\partial}
\newcommand{\Int}{{\rm Int}}
\renewcommand{\d}{\delta}
\newcommand{\e}{\epsilon}
\newcommand{\br}{{\Bbb R}}
\renewcommand{\r}{r}
\newcommand{\ga}{\Gamma_a}
\newcommand{\gb}{\Gamma_b}
\newcommand{\rg}{\hat \Gamma}
\newcommand{\gap}{\Gamma_a^+}
\newcommand{\gbp}{\Gamma_b^+}
\newcommand{\rga}{\hat \Gamma_a}
\newcommand{\rgb}{\hat \Gamma_b}
\newcommand{\rgap}{\hat \Gamma_a^+}
\newcommand{\rgbp}{\hat \Gamma_b^+}
\newcommand{\rgaa}{\hat \Gamma_1}
\newcommand{\rgbb}{\hat \Gamma_2}
\newcommand{\rgaap}{\hat \Gamma_1^+}
\newcommand{\rgbbp}{\hat \Gamma_2^+}
\newcommand{\he}{\hat e}
\newcommand{\na}{n_{a}}
\newcommand{\nb}{n_{b}}
\newcommand{\sign}{\text{\rm sign}}
\newcommand{\BH}{\Bbb H^3}

\input epsf.tex

\section{Introduction}

A Conway sphere for a knot $K$ in $S^3$ is a sphere $S$ which
intersects $K$ at $4$ points, such that punctured sphere $S - K$ is
incompressible in $S^3 - K$.  In this case the sphere $S$ cuts $(S^3,
K)$ into two non-splittable tangles $(B_1, \tau_1)$ and $(B_2,
\tau_2)$, where $B_i$ is a 3-ball, and $t_i$ is a pair of properly
embedded arcs in $B_i$.  An arborescent knot is obtained by gluing
rational tangles together in various ways.  See for example [Wu2] or
[Ga].  An arborescent knot $K$ is {\it large\/} if it has a Conway
sphere.  It is known [HT, Oe, Wu2] that $K$ is large if and only if
its complement is large in the sense that it contains an embedded
closed essential surface.

A nontrivial Dehn surgery on a hyperbolic knot $K$ in $S^3$ is {\it
  exceptional\/} if the resulting manifold is either reducible,
toroidal, or a small Seifert fiber space.  The Geometrization
Conjecture asserts that non-exceptional surgeries yield hyperbolic
manifolds.  Thurston [Th] showed that a hyperbolic knot admits only
finitely many exceptional surgeries.  

All large arborescent knots are hyperbolic.  It is known that most
large arborescent knots admit no exceptional surgery.  Define $T(r,
s)$ to be a Montesinos tangle which is the sum of two rational tangles
associated to rational numbers $r, s$ respectively.  See Section 2 for
more details.  A knot $K$ is an arborescent knot of type II if it has
a Conway sphere cutting it into two Montesinos tangles of type $T(r_i,
\, 1/2)$.  It was shown in [Wu2, Theorem 3.6] that if a large
arborescent knot $K$ is not of type II then all nontrivial surgeries
on $K$ are Haken and hyperbolic, so there is no exceptional surgery on
$K$.  When $K$ is an arborescent knot of type II, all non-integral
surgeries are Haken and hyperbolic, and all integral surgeries are
laminar in the sense that the resulting manifolds contain essential
laminations; in particular, it is irreducible.  It remains to
determine which type II knots admit integer surgeries producing
toroidal or Seifert fibred manifolds.  The following is our main
theorem, which determines all such knots and the exceptional surgeries
on them.  The knot $K_1$ in the theorem is given in Figure 2.2(b).
$K_3$ is actually the mirror image of $K_1$, so there are essentially
only two knots up to homeomorphism of $(S^3,K)$.  $K_2$ is obtained
from $K_1$ by changing crossings on the right half of the diagram of
$K_1$ in Figure 2.2(b).

\begin{thm} 
  Let $K_1, K_2, K_3$ be the three knots in Definition 2.3.  Let $K$
  be a large arborescent knot in $S^3$, and let $\delta$ be a
  non-meridional slope on $\bdd N(K)$.  Then $K(\delta)$ is an
  exceptional surgery if and only if $(K, \delta)$ is isotopic to
  $(K_1, 3)$, $(K_2, 0)$ or $(K_3, -3)$, in which case $K(\delta)$ is
  toroidal.
\end{thm}

The paper is organized as follows.  In Section 2 we will define some
special disks in rational tangle spaces.  These are the pieces which
will be used to build the toroidal surfaces in the exteriors of the
knots in Theorem 1.1.  Section 3 defines an index $i(G, Q)$ for a
surface $G$ relative to $Q$, and proves its additivity and some other
properties.  Now let $K$ be a type II knot, which is the union of two
Montesinos tangles $T_i = T(r_i, 1/2)$.  Let $F$ be a punctured
essential torus in the exterior of $K$ with integer boundary slope,
and let $F_i$ be the intersection of $F$ with the tangle space of
$T_i$.  The important fact is that $F$ can be chosen so that each
component $G$ of $F_i$ must have zero index relative to the tube
around the unknotted string of $T_i$.  It will then be shown in
Section 4 that $G$ must be a special surface in the sense that it is
the union of special disks in the two rational tangle spaces.  This
quickly leads to a proof that $r_i \equiv \pm 1/3$ mod $1$.  In
Section 5 we define the relative framing for a surface in a tangle
space.  They can be used to calculate the boundary slope of the
surface $F = F_1 \cup F_2$, and show that if $F$ has integer boundary
slope then its intersection with the punctured Conway sphere must have
a very special configuration, which completely determines the gluing
map between the two tangles.  It will follow that if $K(\delta)$ is
toroidal then $K$ must be one of the three knots in the theorem, and
there is only one possible choice of $F$.  In Section 6 it will be
shown that the surface $F$ constructed is incompressible and
$\bdd$-incompressible, and it gives rise to an essential torus in the
surgered manifold.  This, together with some known results about
surgery on type II knots, will complete the proof of Theorem 1.1.

Unless otherwise stated, a surface in a 3-manifold $M$ is assumed to
be either on the boundary of $M$ or properly embedded in $M$.  Given a
set $X$ in a surface or 3-manifold, denote by $N(X)$ a regular
neighborhood of $X$, and by $|X|$ the number of components in $X$.  If
$P$ is a surface in a 3-manifold $M$, denote by $M|P$ the manifold
obtained by cutting $M$ along $P$.

\section{Special disks in tangle spaces}

A {\it tangle\/} $T$ is a pair $(B, \tau)$, where $B$ is a 3-ball with
four specified points on $\bdd B$, and $\tau = \tau_1 \cup \tau_2$ is
a pair of arcs in $B$ connecting these points.  We will identify $B$
with either a pillow case with $\bdd \tau$ the four corners, or the
one point compactification of the lower half space, i.e.\ $B = \hat
R^3_- = \{(x,y,z)\, | \, z\leq 0\} \cup \{\infty\}$, with $\bdd \tau$
identified to the four points $(\pm 1, \pm 1)$ and the front side of
the pillow case identified to the square $[-1,1] \times [-1, 1]$ in
$\hat R^2 = R^2 \cup \{\infty\}$, which is identified to $(R^2
\times \{0\}) \cup \{\infty\} = \bdd \hat R^3_-$.

Given a tangle $T = (B, \tau)$, let $E(T)$ be the closure of $B -
N(\tau)$, called the tangle space of $T$.  Let $S(T) = \bdd B = \hat
R^2$, and let $P(T)$ be the 4-punctured sphere $\bdd B \cap E(T)$.
Let $S_+(T)$ (resp.\ $S_-(T)$ be the closure of the right (resp.\
left) half plane of $\hat R^2 = \bdd B$.  Similarly, define
$P_{\pm}(T) = S_{\pm}(T) \cap E(T)$, which is a twice punctured disk.
Denote by $U(T) = U_+(T) \cup U_-(T)$ the two annuli $\bdd E(T) - \Int
P(T)$, with $U_+(T)$ the one containing the upper right component of
$\bdd P$ on $\hat R^2$.  We now have a decomposition of the boundary
of $E(T)$ as follows.  $$\bdd E(T) = P_-(T) \cup P_+(T) \cup U_-(T)
\cup U_+(T).$$

We refer the readers to [Co] or [Wu1] for the definition of rational
tangles.  Roughly speaking, the strings $\tau$ of a rational tangle of
slope $r=p/q$ is obtained by pushing into the interior of $B$ the
interior of two arcs of slope $r$ on the boundary of the pillow case
$B$ connecting the four corners of $B$.  {\it Throughout this paper we
will always assume that $q \geq 2$.}

Given two tangles $T_i = (B_i, \tau_i)$, we may construct a new tangle
$T_1 + T_2$ by identifying the disk $S_+(T_1) \subset \hat R^2$ with
$S_-(T_2) \subset \hat R^2$ using the map $(x,y) \to (-x, y)$ and then
identifying $B_1 \cup B_2$ to $B = \hat R^3_-$ so that a boundary
point of $B_1$ or $B_2$ on $\bdd (B_1 \cup B_2)$ is mapped to the
point with the same coordinates on $\bdd B = \hat R^2$.  Denote by
$T(r_1, r_2)$ the Montesinos tangle $T(r_1) + T(r_2)$.

Two tangles $T_1, T_2$ are {\it weakly equivalent\/} if $T_1$ can be
deformed to $T_2$ by an isotopy $\varphi_t$ of $B$.  They are {\it
  $P$-equivalent\/} if the isotopy $\varphi_t$ above is rel $\bdd
P_+(T_1)$, and {\it equivalent\/} if $\varphi_t$ is rel $S(T)$.  Thus
for example, $T(r)$ is $P$-equivalent to $T(r+k)$ for any integer $k$,
and $T(1/3, -1/2)$ is $P$-equivalent to $T(1/3, 1/2)$.  Two
tangles are considered the same if they are equivalent.  

A surface $F$ in $E(T)$ is {\it tight\/} if $\bdd F \neq \emptyset$,
and it intersects each of $P_{\pm}$ and $U_{\pm}$ in essential arcs or
essential circles.  In this case $\bdd F$ is a set of essential loops
on $\bdd E(T)$.  Thus a tight disk is an essential disk in $E(T)$
(i.e.\ a compressing disk of $\bdd E(T)$), and any essential disk in
$E(T)$ is isotopic to a tight disk.

\bigskip
\leavevmode

\centerline{\epsfbox{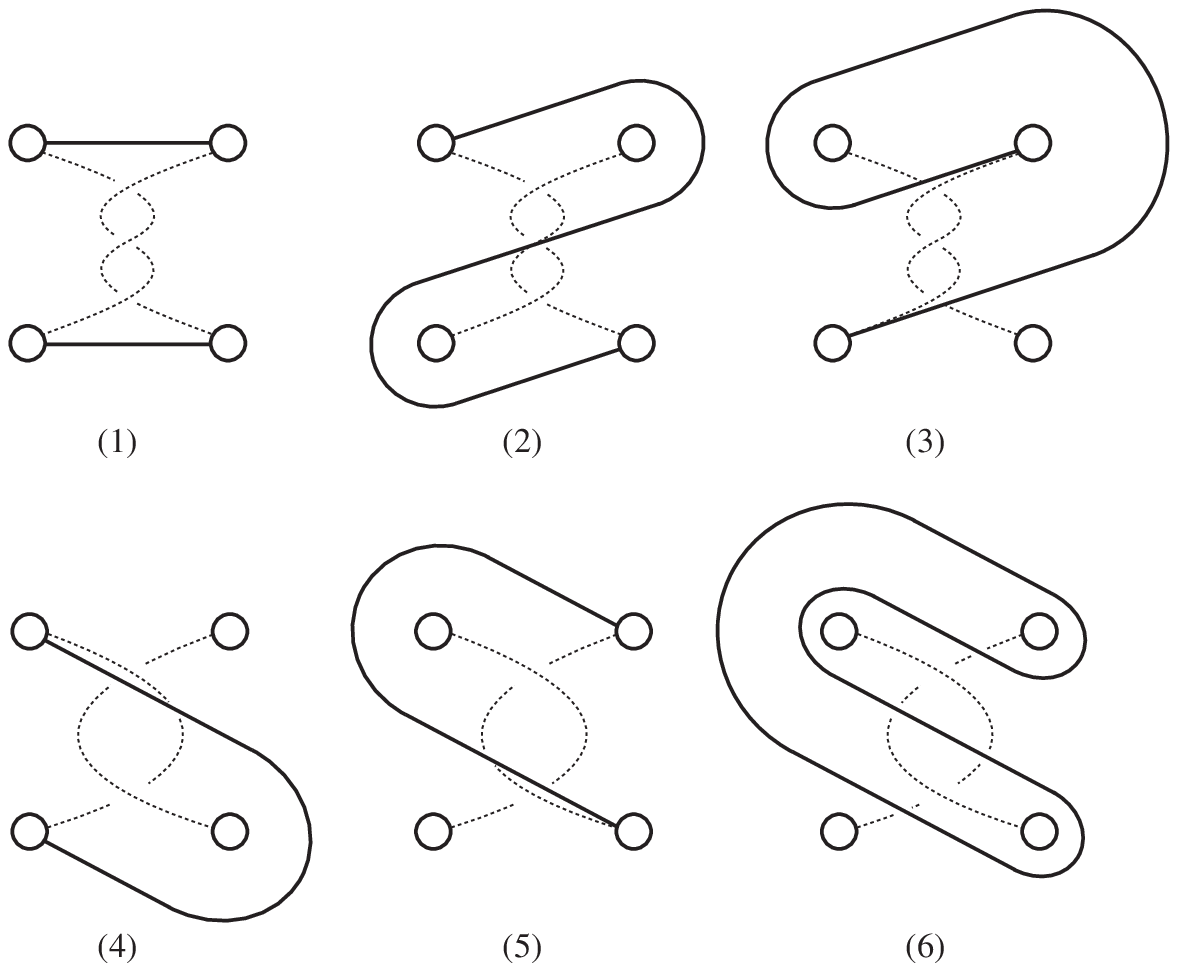}}
\bigskip
\centerline{Figure 2.1}
\bigskip

The strings $\tau$ in a rational tangle $T(p/q)=(B, \tau)$ are rel
$\bdd \tau$ isotopic to a pair of arcs $\tau' = \tau'_1 \cup \tau'_2$
on $\bdd B$.  To be specific, let $\tau'_2$ be the one with an
endpoint at $(1,1) \in \hat R^2$.  Let $\tau'_0$ be a pair of
horizontal arcs on $R^2$ connecting $\bdd \tau$.  Let $c_i = \tau'_i
\cap P(T)$ for $i=0,1,2$, and let $c_3$ be the curve on $P$ that
separates $c_1$ and $c_2$, which is unique up to isotopy.  Thus for
example, for $T(1/3)$ the curves $c_0, c_1, c_2$ are shown in Figure
2.1(1)--(3), and for $T(-1/2)$ the curves $c_1, c_2, c_3$ are shown in
Figure 2.1(4)--(6), respectively.  A disk $D$ in $E(T)$ is a special
disk if its intersection with $P(T)$ is one of these curves.  It is
further required that $q\leq 3$ if $D\cap P(T) = c_1, c_2$, and $q=2$
if $D\cap P(T) = c_3$.  More explicitly, we have the following
definition.  Note that the curve in Figure 2.1(k) is the boundary
curve of a type (k) special disk.

\begin{defn} {\rm Let $T = T(\pm 1/q)$.  A disk $D$ in $E(T)$ is a
  {\it special disk\/} if it is one of the following type.

Type (1).  $T = T(\pm 1/q)$, $q$ odd, and $D \cap P(T) = c_0$.  

Type (2).  $T = T(\pm 1/3)$, and $D \cap P(T) = c_1$.  

Type (3).  $T = T(\pm 1/3)$, and $D \cap P(T) = c_2$.  

Type (4).  $T = T(\pm 1/2)$, and $D \cap P(T) = c_1$.  

Type (5).  $T = T(\pm 1/2)$, and $D \cap P(T) = c_2$.  

Type (6).  $T = T(\pm 1/2)$, and $D \cap P(T) = c_3$.  
 }
\end{defn}

\begin{lemma} Let $T = T(p/q)$ with $q = 2$ or odd.  Let $Q = P_+(T)$
  if $q\geq 3$, and $Q = P_-(T) \cup U_+(T)$ if $q=2$.  Suppose $D$ is
  a tight disk in $T(p/q)$ such that $|D \cap Q| \leq 2$.  Then $T$ is
  $P$-equivalent to $T(\pm 1/q)$, and $D$ is a special disk.  In
  particular, $Q$ and $Q' = \bdd E(T) - \Int Q$ are incompressible,
  and there is no disk in $E(T)$ intersecting each of $Q$ and $Q'$ at
  a single essential arc.
\end{lemma}

\proof This is essentially [Wu1, Lemmas 2.1 and 2.2].  Let $P =
P_+(T)$ if $q\geq 3$, and $P = P_-(T)$ if $q=2$.  Clearly $D$ must
intersect $P$ in a nonempty set of arcs as otherwise we would have $T
= T(1/0)$.  Since $|D \cap P| \leq 2$, $D$ is a monogon or bigon as
defined in [Wu1], which have been classified in Lemmas 2.1 and 2.2 of
[Wu1].  When it is a monogon $D$ is a special disk of type (5).
Bigons appear when $T$ is a torus tangle or wrapping tangle (see [Wu1]
for definition), or a twist tangle; but since $T$ is rational and
$q=2$ or odd, the first two cases does not happen.  Thus from the
proof of [Wu1, Lemma 2.2] we see that if $D$ is a bigon then it is one
of the types (1), (2), (3), (4) or (6) in Definition 2.1, or $T$ is
$P$-equivalent to $T(1/4)$ and $D$ intersects $P$ in two arcs with
boundary on the outer component of $\bdd P$, but since we have assumed
that $q=2$ or odd, the last case is impossible.  \qed

\bigskip
\leavevmode

\centerline{\epsfbox{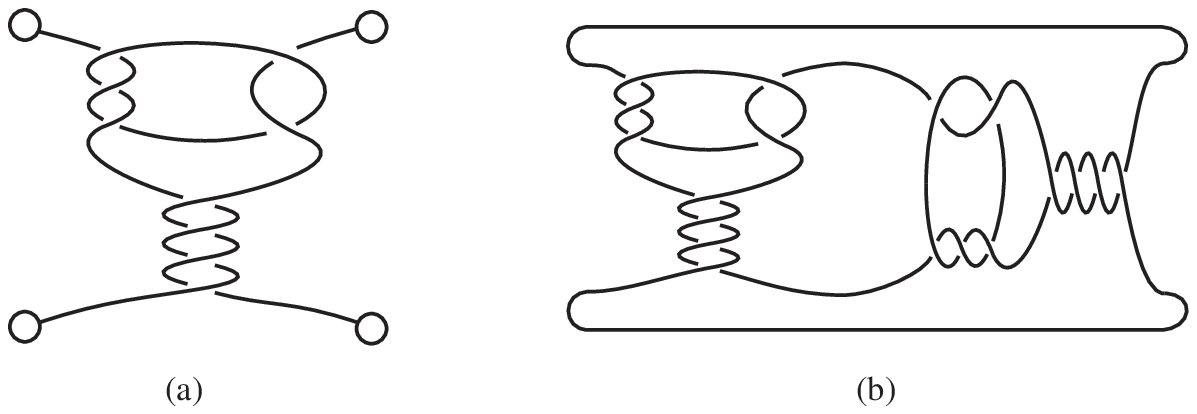}}
\bigskip
\centerline{Figure 2.2}
\bigskip

Denote by $T(r_1, r_2;\, n)$ the tangle obtained from $T(r_1, r_2) =
(B_3, \tau)$ by twisting the two lower endpoints of $\tau$ by $n$
left hand half twists.  See Figure 2.2(a) for the tangle $T(1/3,
-1/2;\, 4)$.

\begin{defn} {\rm
Let $\eta: \hat R^2 \to \hat R^2$ be the map which is a $\pi/2$
counterclockwise rotation about the origin followed by a reflection
along the $y$-axis.  Define three knots $K_1, K_2, K_3$ by

$(S^3, K_1) = T(1/3, -1/2;\, 4) \cup_{\eta} T(1/3, -1/2;\, 4)$, 

$(S^3, K_2) = T(1/3, -1/2;\, 4) \cup_{\eta} T(-1/3, 1/2; -4)$, and 

$(S^3, K_3) = T(-1/3, 1/2; -4) \cup_{\eta} T(-1/3, 1/2; -4)$.  }
\end{defn}

Alternatively, $K_i$ can be obtained by shifting the first tangle to
the left, the second tangle to the right, rotating the second tangle
counterclockwise by an angle of degree $\pi/2$, and then connecting
the endpoints of the tangles by arcs on $R^2$ which are horizontal
except near the endpoints of the strings of the tangles, as shown in
Figure 2.2(b) for $K_1$.  Note that $K_3$ is the mirror image of
$K_1$, and $K_2$ is obtained from $K_1$ by taking mirror image the
right half of $K_1$.  Theorem 1.1 shows that these are the only large
arborescent knots which admit exceptional surgery, and each of them
admits exactly one such surgery.

\section{Index of essential surfaces}

Let $Q$ and $F$ be surfaces in $M$, intersecting in general position.
Denote by $a(F, Q)$ the number of arc components of $F\cap Q$.  The
{\it index of $F$} in $M$ relative to $Q$ is defined as $$i(F, Q) =
\chi(F) - \frac 12 a(F, Q)$$
where $\chi(F)$ is the Euler
characteristic of $F$.  This is the same as the cusped Euler
characteristic defined in [Wu3] for sutured manifold, only that now
$Q$ is not required to be ``cusps'', which by definition is a set of
annuli and tori on the boundary of a 3-manifold.

\begin{lemma} Let $Q$ be a surface on $\bdd M$, and $Q'$ an essential
  surface properly embedded in $M$ and disjoint from $Q$.  Let $M' =
  M|Q'$, and let $Q'_1, Q'_2$ be the two copies of $Q'$ on $\bdd M'$.
  Let $F$ be a surface in $M$, and let $F' = F|Q'$ be the
  corresponding surface in $M'$.  Then $$
  i(F', Q'_1 \cup Q'_2 \cup Q)
  = i(F, Q).$$
\end{lemma}

\proof Put $k = a(F, Q')$.  Note that $Q'_1, Q'_2, Q$ are mutually
disjoint compact surfaces on $\bdd M'$, hence $a(F', Q'_1 \cup Q'_2
\cup Q) = 2k + a(F, Q)$.  Since $Q'$ intersects $F$ in $k$ arcs and
possibly some circle components, after cutting $F$ along $F\cap Q'$ we
have $\chi(F') = \chi(F) + k$.  It follows that
\begin{eqnarray*}
i(F, Q) & = & \chi(F) - \frac 12 a(F, Q) = \chi(F') - k - \frac 12
a(F, Q)  \\
& = &\chi (F') - \frac 12 a(F', Q'_1 \cup Q'_2 \cup Q) \\
& = & i(F', Q'_1 \cup Q'_2 \cup Q). \quad \Box
\end{eqnarray*}

The lemma shows that index is invariant when cutting along a surface
disjoint from $Q$.  We remark that it is important to assume that $Q$
is a compact surface and is disjoint form $Q'$, otherwise the lemma may
not be true.

The most useful case is when $Q'$ is a separating surface, cutting $M$
into $M_1$ and $M_2$.  The following additivity lemma follows
immediately from Lemma 3.1.

\begin{lemma} {\rm (Additivity of index)} 
  Suppose $Q$ is a compact subsurface of $\bdd M$, and $Q'$ is a
  separating surface in $M$ disjoint from $Q$, cutting $M$ into $M_1$
  and $M_2$.  Let $Q'_i$ be the copy of $Q'$ on $\bdd M_i$, $Q_i = (Q
  \cap M_i) \cup Q'_i$, and $F_i = F \cap M_i$.  Then $$
  i(F, Q) =  i(F_1, Q_1) + i(F_2, Q_2). \quad \Box $$
\end{lemma}

\begin{lemma} Let $T = T(p/q)$ with $q=2$ or odd.  Let $Q = P_+(T)$
  for $q>2$, and $Q = P_-(T) \cup U_+(T)$ for $q=2$.  Let $F$ be a
  tight disk in $E(T)$.  Then $i(F, Q) \leq 0$, and equality holds if
  and only if $T$ is $P$-equivalent to $T(\pm 1/q)$ and $F$ is a
  special disk.
\end{lemma}

\proof It is easy to check that each special disk in Definition 2.1
intersects $Q$ in two arcs and hence has $i(F, Q) = 0$.  Conversely,
if $i(F, Q) \geq 0$ then $|F \cap Q|\leq 2$, so by Lemma 2.2 it is a
special disk and hence $i(F,Q)=0$.  
\qed

\begin{lemma} Suppose $T=T(p/q, 1/2)$ and $F$ is a tight surface in
  $E(T)$.  Let $U = U_+(T)$.   Then $i(F, U) \leq 0$.  In particular,
  $P(T)$ and $\bdd E(T) - U$ are incompressible, and there is no disk
  in $E(T)$ intersecting $U_+(T)$ at a single essential arc.
\end{lemma}

\proof If $i(F,U)>0$ then $F$ is a disk and $F\cap U$ has at most one
arc component.  Isotope the decomposition surface $P$ in $E(T) =
E(T(p/q)) \cup _P E(T(1/2))$ so that $F\cap P$ is minimal.  Then an
innermost circle outermost arc argument would lead to a contradiction
to Lemma 3.3 because one of the (at least two) disks cut off by
outermost arcs will be disjoint from $U$ and hence is isotopic to a
tight disk with positive index.
\qed

\section{Special surfaces in knot exterior}

Throughout this section we assume that $K$ is a type II knot which is
the union of two length 2 Montesinos tangles $T_i = (B_i, \tau_i)=
T(p_i/q_i, 1/2) = T_{i1} + T_{i2}$, where $T_{i1} = T(p_i/q_i)$ and
$T_{i2} = T(1/2)$.

Let $S = \bdd B_i$, and let $P = P(T_1) = P(T_2) = S \cap E(K)$, which
cuts $E(K)$ into $E(T_1)$ and $E(T_2)$.  Let $P_i = E(T_{i1}) \cap
E(T_{i2})$ be the twice punctured disk cutting $E(T_i)$ into
$E(T_{i1})$ and $E(T_{i2})$.  Write $T_{ij} = (B_{ij}, \tau_{ij})$.
Thus $E(K)$ is the union of four rational tangle spaces $E(T_{ij})$,
$i,j=1,2$.  Let $U_i = U_+(T_i)$, which is the component of $U(T_i)$
lying in $E(T_{i2})$.

\begin{defn} {\rm A surface $F$ in $E(K)$ or $E(T_i)$ is a {\it
      special surface\/} if it intersects each $E(T_{ij})$ in special
    disks.}
\end{defn}

The purpose of this section is to show that if $F$ is an essential
punctured torus in $E(K)$ with integer slope then it is a special
surface up to isotopy.  We will then use this result to show that $q_i
= 3$ for $i=1,2$.

If $F$ is a compact surface in $E(T)$ intersecting $U(T)$ in arcs on
$\bdd F$, then $F - U(T)$ is $F$ with the arcs $F\cap U(T) \subset
\bdd F$ removed, so it is non-compact.  An arc $\beta$ in $F-U(T)$ is
considered to be {\it essential\/} if it does not cut off a compact
disk from $F- U(T)$.  Thus for example, if $F$ is a disk intersecting
$U(T)$ in two arcs on $\bdd F$ then there is exactly one essential arc
on $F-U(T)$ up to isotopy.

A {\it $P$-compressing disk\/} of a surface $F$ in $E(T)$ is a disk
$D$ in $E(T)$ such that $\bdd D = \alpha \cup \beta$, where $\alpha$
is an arc on $P$ and $\beta = D \cap (F-U(T))$ is an arc on $F - U(T)$
which is essential in the above sense.  If such a disk exists then $F$
is {\it $P$-compressible}, otherwise it is {\it $P$-incompressible}.

\begin{lemma} Let $F$ be a punctured essential torus in $E(K)$ such
  that $\bdd F$ has integer slope on $\bdd E(K)$, and the complexity
  $(|F \cap \bdd P|, |F\cap P|)$ is minimal in lexicographic order.
  Let $F_j$ be a component of $F\cap E(T_i)$.  Then (i) $F_j$ is
  tight, (ii) $i(F_j, U_i) = 0$, and (iii) $F_j$ is incompressible and
  $P$-incompressible.
\end{lemma}

\proof (i) Since $|F \cap \bdd P|$ is minimal, $F$ intersects each of
the four annuli $\bdd N(K) | \bdd P$ in essential arcs.  By [Wu1,
Lemma 3.3] each $E(T_i)$ is a handlebody and hence irreducible, so if
$F\cap P$ contains a trivial loop on $P$ then an innermost circle
argument would show that one could isotope $F$ to reduce $|F \cap P|$
without increasing $|F \cap \bdd P|$, which is a contradiction.  If
$F\cap P$ has a trivial arc on $P$ then an outermost one would cut off
a $\bdd$-compressing disk for $F$, which is impossible because $F$ is
essential.

(ii) Since $P$ is incompressible [Wu1, Lemma 3.3], each circle
component of $F\cap P$ is also essential on $F$.  It follows that each
boundary component of $F\cap E(T_i)$ is a nontrivial loop on $\bdd
E(T_i)$, so each disk component of $F \cap E(T_i)$ is an essential
disk in $E(T_i)$.

Let $\hat F$ be the torus obtained from $F$ by capping off each
boundary component of $F$ with a disk.  Define a graph $\Gamma$ on
$\hat F$ with the attached disks as fat vertices and $P\cap F$ as
edges.  A component of $P \cap F$ is either a {\it circle edge\/} or
an {\it arc edge\/} of $\Gamma$, depending on whether it is a circle
or arc.  Note that a circle edge is not incident to any vertex, so
$\Gamma$ is actually a graph in which some edges are loops without
vertices.  

Since $\bdd F$ has integer boundary slope and $P$ has $4$ boundary
components, each vertex $v$ of $\Gamma$ has valence $4$.  Denote by $C$
the number of corners at all vertices of $\Gamma$ which lie in either
$U_1$ or $U_2$, where $U_i = U_+(T_i)$, which is the component of
$\bdd N(K) \cap E(T_i)$ lying in the $E(T_{i2}) = E(T(1/2))$ part.
Denote by $V$ and $E$ the numbers of vertices and {\it arc\/} edges of
$\Gamma$, respectively.  Since the boundary of each vertex of $\Gamma$
travels through each of $U_1$ and $U_2$ once, we have $C = 2V$.  Since
each vertex is incident to $4$ arc edge endpoints, we also have $E =
2V$.  Let $F_j$ be the faces of $\Gamma$, and assume it lies in
$E(T_i)$.  The Eular characteristic formula, which one can check is
still valid for graphs with circle edges, gives
\begin{eqnarray*}
0 & = & V - E + \sum \chi(F_j) = - V + \sum \chi(F_j) = - \frac 12 C +
\sum \chi(F_j) \\
& = & \sum (\chi(F_j) - \frac 12 |F_j \cap (U_i)|) = \sum
i(F_j, U_i).
\end{eqnarray*}

By (i) and Lemma 3.4 we have $i(F_j, U_i) \leq 0$ for each component
$F_j$ of $F\cap E(T_i)$.  Since these are exactly the faces of
$\Gamma$, it follows from the above that $i(F_j, U_i) = 0$ for all
$F_j$.

(iii) Since $F$ is incompressible and $|F \cap P|$ is minimal, it is
easy to see that $F_j$ is also incompressible.  By (ii) $F_j$ is
either a disk intersecting $U_i$ twice or an annulus disjoint from
$U_i$.  Suppose $D$ is a $P$-compressing disk for $F_j$ as in the
definition, and let $\alpha = D \cap P$.  Since $\beta = |D \cap F_j|$
is essential on $F_j - U_i$, the two points $\bdd \alpha = \bdd \beta$
lie on distinct components of $F_j \cap P$, hence an isotopy of $F$
via $D$ would create a surface $F'$ which has the same or smaller
complexity than $F$ and yet one of the faces $F'_j$ deformed from
$F_j$ has $i(F'_j, U_i)>0$, which is a contradiction to (ii).  \qed

\begin{lemma} Let $F$ be an essential punctured torus in $E(K)$ with
  integer boundary slope.  Then $F, P_1, P_2$ can be isotoped so that
  $F$ is a special surface.
\end{lemma}

\proof We may assume that $(|F\cap \bdd P|, |F \cap P|)$ is minimal up
to isotopy, so Lemma 4.2 applies, and we have $i(F_k, U_i)=0$ for any
component $F_k$ of $F \cap E(T_i)$.  We will show below that $P_i$ can
be isotoped so that each component $D$ of $F\cap E(T_{ij})$ is a tight
disk.  In that case by Lemma 3.3 we will have $i(D, U_i) \leq 0$;
since $F_k$ is a union of such disks, by the Additivity Lemma 3.2 we
will have $i(D, U_i)=0$ for all $D$, and hence by Lemma 3.3 $D$ is a
special disk in $E(T_{ij})$, which will complete the proof.

Let $Q_i = F \cap E(T_i)$.  Isotope $P_i$ so that $(|Q_i \cap \bdd
P_i|, |Q_i \cap P_i|)$ is minimal.  This implies that for $R =
P_-(T_{i1})$, $P_+(T_{i2})$ or a component of $U(T_{ij})$, each arc
component of $Q_i \cap R$ is essential on $R$.  By Lemma 4.2 $Q_i$ is
incompressible and $P$-incompressible, so the above minimality and the
$\bdd$-incompressibility of $F$ imply that each component of $Q_i \cap
P_i$ is also essential on $P_i = P_+(T_{i1}) = P_-(T_{i2})$.
Therefore all components of $F\cap E(T_{ij})$ are tight.  It remains
to show that each component of $F \cap E(T_{ij})$ is a disk.

By Lemma 2.2 $P_i$ is incompressible and $\bdd$-incompressible in
$E(T_i)-U_i$, so no component of $P_i\cap Q_i$ is a trivial loop on
$Q_i$, or an arc which is trivial on $Q_i - U_i$ in the sense that it
cuts off a disk on $Q_i$ disjoint from $U_i$.  Thus if a component
$F_k$ of $Q_i$ is a disk then $P_i$ cuts $F_k$ into disks.  Similarly
if $F_k$ is an annulus then no component of $F_k \cap P_i$ is an
inessential arc or inessential circle on $F_k$.  It now suffices to
show that $\bdd F_k \cap P_i \neq \emptyset$ because then $F_k \cap
P_i$ is a nonempty set of essential arcs on $F_k$, cutting it into
disks.

Let $F_k$ be an annulus component of $Q_i$ with $\bdd F_k$ disjoint
from $P_i$.  This implies that $F_k \cap U_-(T) = \emptyset$ because
each component of $F_k \cap U_-(T_i)$ is a component of $F\cap
U_-(T_i)$, which intersects $P_i$ at two points.  Since $i(F_k,
U_i)=0$, we also have $F_k \cap U_+(T_i) = \emptyset$.  Therefore
$\bdd F_k \subset P$.  Note that no component $C'$ of $\bdd F_k$ is
parallel to a component $C$ of $\bdd P$ as otherwise the arc
components of $F\cap P$ with endpoints on $C$ would lie in the annulus
between $C$ and $C'$ and hence would be trivial arcs, and there are
$|\bdd F|/2 > 0$ of them, which can be used to $\bdd$-compress $F$, a
contradiction to the fact that $F$ is essential.  Therefore the two
components of $\bdd F_k$ are both parallel to the circle $P_i \cap P$
and hence bound an annulus $A$ on $P$.  By Lemma 4.2(iii) $F_k$ is
incompressible, and by [Wu1, Lemma 3.3] $E(T_i)$ is a handlebody, so
$F_k \cup A$ must bound a solid torus $V$ in $E(T_i)$.  Now $A$ cannot
be meridional on $V$ because $P$ is incompressible, and $A$ cannot run
more than once along the longitude of $V$ as otherwise the union of
$V$ and a regular neighborhood of a disk on $\bdd B_i$ bounded by a
component of $\bdd F_k$ would be a punctured lens space in the 3-ball
$B_i$, which is absurd.  It follows that $A$ is longitudinal on $\bdd
V$ and hence $F_k$ is $P$-compressible, which contradicts Lemma
4.2(iii).  \qed

\begin{prop} Suppose $K$ is an arborescent knot of type II, and
  suppose $E(K)$ contains an essential punctured torus with integer
  boundary slope.  Then $K$ is the union of $T_1, T_2$, each of which
  is weakly equivalent to $T(1/3, 1/2)$ or $T(-1/3, 1/2)$.
\end{prop}

\proof By definition $K$ is the union of $T_1, T_2$ with $T_i =
T(p_i/q_i, 1/2) = T_{i1} + T_{i2}$.  By Lemma 4.3 we may assume that
$F$ intersects each $E(T_{ij})$ in special disks.  If some $q_i>3$
then the only special disk in $E(T_{i1})$ is of type (1), hence each
component of $F \cap P_i$ is an arc with only one endpoint on the
circle $P_i \cap P$, where $P_i$ is the twice punctured disk
$E(T_{i1}) \cap E(T_{i2})$.  From the definition we see that the only
special disks in $E(T_{i2}) = E(T(1/2))$ intersecting $P_i$ in such
arcs are those of type (5), which are disjoint from $U_+(T_2)$.
Therefore $F \cap U_+(T_i) = \emptyset$, which is a contradiction
because $F$ must intersect each $U_{\pm}(T_i)$ in exactly $|\bdd F|>0$
arcs.  \qed

\section{Boundary slopes of special surfaces}

Let $T=(B^3, \tau)$ be a tangle, with $B^3$ identified to $\hat
R^3_-$.  We assume that $\tau$ is oriented.  Let $\alpha$ be a pair of
arcs on $R^2 - \Int I^2$, shown in Figure 5.1(1) when the orientations
of $\tau$ are opposite near the two upper endpoints, or in Figure
5.1(2) otherwise.  Then $\hat \tau = \tau \cup \alpha$ is a link in
$R^3$ with orientation induced from that of $\tau$.

\bigskip
\leavevmode

\centerline{\epsfbox{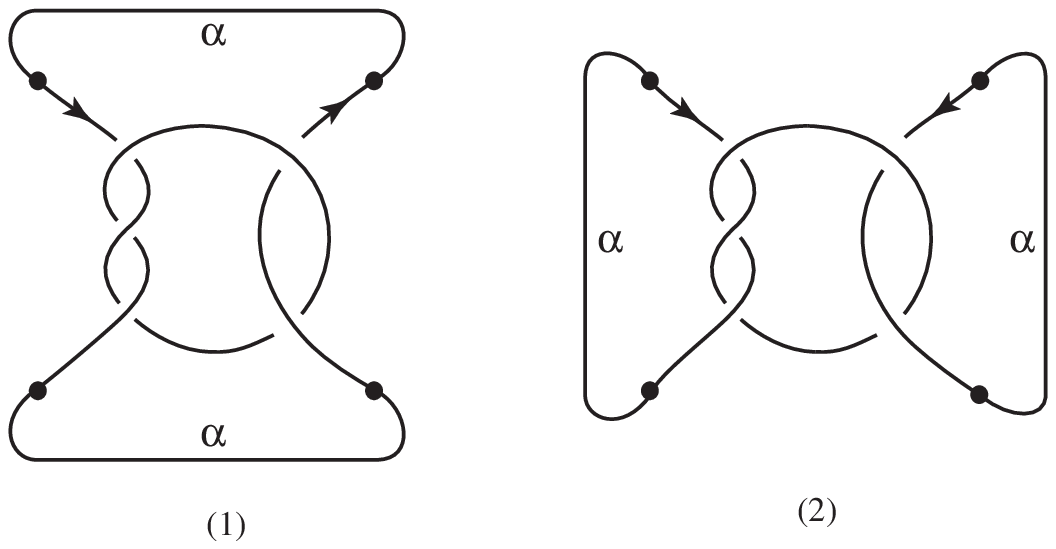}}
\bigskip
\centerline{Figure 5.1}
\bigskip

Let $\gamma$ be a set of essential arcs on $U(T)$, and assume that $n
= |\gamma \cap U_-(T)| = |\gamma \cap U_+(T)|$.  Each component of
$\gamma$ is isotopic in $N(\tau)$ to a component of $\tau$, so the
orientation of $\tau$ induces an orientation on $\gamma$.  Let $p:
R^3_- \to R^2 = \bdd R^3$ be the standard projection.  We always
assume that $\tau$ and $\gamma$ are in {\it regular position\/} in the
sense that

(i) $p(\tau) \subset I^2$, where $I = [-1, 1]$; 

(ii) $p: \tau\cup \gamma \to R^2$ is an immersion with only double
crossings; and 

(iii) $p(\gamma) \cap \alpha = \emptyset$.

Condition (iii) is to guarantee that all crossings between $p(\gamma)$
and $p(\tau \cup \alpha)$ appear inside of the square $I^2$, so if we
close up $\gamma$ by arcs parallel to $\alpha$ to obtain $\hat
\gamma$ then the linking number between $\hat \gamma$ and $\hat \tau =
\tau \cup \alpha$ can be calculated using crossings between $\tau$ and
$\gamma$.

Since $\tau$ and $\gamma$ are oriented, each crossing is assigned a
sign according to the right hand rule, as given in Rolfsen's book
[Ro, p.132].

\begin{defn} {\rm (1) The {\it relative linking number\/} between
    $\tau$ and  $\gamma$, denoted by $lk(\tau, \gamma)$, is the sum of
    the signs of crossings at which $\gamma$ passes below $\tau$. 

    (2) Let $F$ be a surface in $E(T)$ such that $F \cap U(T)$ is in
    regular position.  Then the {\it relative framing\/} of $F$ in
    $E(T)$ is defined as $\theta(F) = lk(\tau, F\cap U(T))$.  }
\end{defn}

An isotopy of $\tau \cup \gamma$ is a {\it regular isotopy\/} if $\tau
\cup \gamma$ is in regular position at any time during the isotopy.
Similarly an isotopy of a surface $F$ in $E(T)$ is a regular isotopy
if its restriction to $F\cap U(T)$ is a regular isotopy.

Consider the four disks $\cup D_i = R^2 - \Int P(T)$.  Since $\gamma$
has $n$ endpoints on each $\bdd D_i$, we can connect $\bdd \gamma$ by
$2n$ arcs $\gamma'$ on $\bdd N(\hat \tau)$ which lies in the upper
half space $R^3_+$.  Define $\hat \gamma = \gamma \cup \gamma'$, with
orientation induced by that of $\gamma$.  Since $\hat \tau \cup \hat
\gamma$ is an oriented link, the linking number $lk(\hat \tau, \hat
\gamma)$ is well defined.

\begin{lemma}  (1) $lk(\tau, \gamma) = lk(\hat \tau, \hat \gamma)$.

  (2) $lk(\tau, \gamma)$ and $\theta(F)$ are regular isotopy
  invariants .

  (3) Let $\psi$ be a rotation of $R^3_-$ along the $z$-axis by an
  angle of $\pi/2$, deforming $\tau$ to $\tau'$ and a surface $F$ to
  $F'$.  Then $\theta(F) = \theta(F')$.

  (4) Suppose $F_j$ are the components of $F$ in $E(T)$.  Then
  $\theta(F) = \sum \theta(F_j)$.
\end{lemma}

\proof (1) It is well known that $lk(\hat \tau, \hat \gamma)$ can be
calculated as the sum of the signs of crossings at which $\hat \gamma$
passes below $\hat \tau$.  By definition $\gamma$ does not pass below
$\alpha$, and $\gamma'$ does not pass below $\hat \gamma$ because it
lies in $R^3_+$ while $\hat \gamma$ lies in $R^3_-$.  Therefore the
crossings at which $\hat \gamma$ passes below $\hat \tau$ are exactly
where $\gamma$ passes below $\tau$, and the result follows.

(2) A regular isotopy does not change the relative position of $\bdd
\tau$ to $\alpha$, hence it extends to an isotopy of $\hat \tau \cup
\hat \gamma$, and the result follows from (1) because the linking
number of a link is an isotopy invariant.

(3) This follows from the definition, because $\psi$ gives a sign
preserving one to one correspondence between the crossings.  

(4) This follows from definition.
\qed

Let $T=T_1 + T_2$, where $T_1 = T(1/3)$ and $T_2 = T(-1/2)$.  Let $F$
be a special surface in $E(T)$ intersecting each $U_{\pm}(T)$ in $n$
arcs.  By definition each component of $F \cap E(T_i)$ is a special
disk.  We may assume that $F$ is a union of $a_i$ copies of $A_i$ for
$i=1,...,6$.  From Figure 2.1 we see that $F \cap P(T_i)$ is as shown
in Figure 5.2($i$), $i=1,2$, where an arc with label $a_j$ indicates
$a_j$ parallel copies of that arc.

\bigskip
\leavevmode

\centerline{\epsfbox{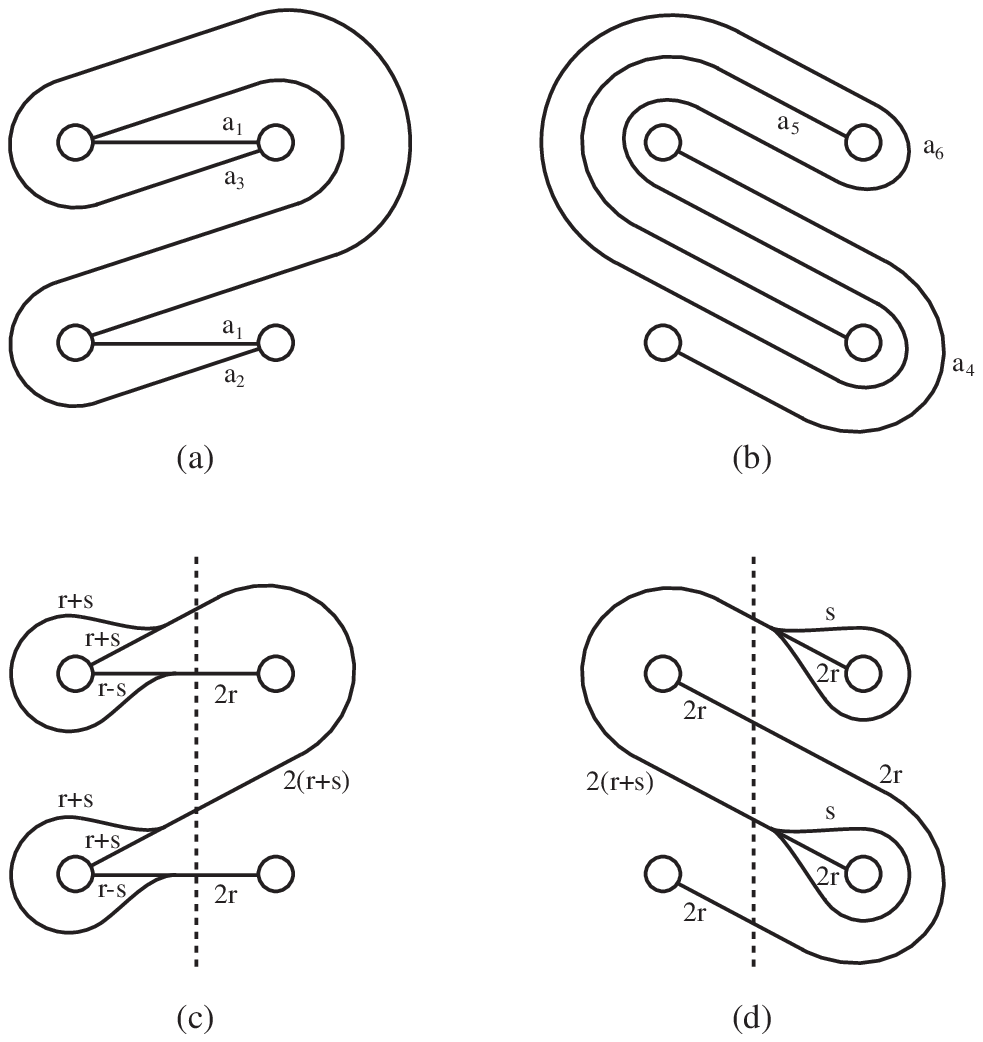}}
\bigskip
\centerline{Figure 5.2}
\bigskip

\begin{lemma} Let $r = n/2$, and $s = a_6$.  Then there is a a special
  surface $F$ in $E(T)$ consisting of $a_k$ copies of disks of type
  (k) and intersecting each $U_i(T)$ at $n$ arcs, if and only if $a_k$
  is a set of non-negative numbers satisfy the following equations.

(1) $n = 2r$ is even;

(2) $a_1 = r - s$;

(3) $a_2 = a_3 = r + s$;

(4) $a_4 = a_5 = a_1 + a_2 = a_1 + a_3 = n$;

(5) $a_6 = s$.
\end{lemma}

\proof Suppose $F$ is such a surface.  Then (5) is from definition,
and (4) follows from the fact that $F$ intersects each boundary
component of $P$ at $n$ points.  By (4) we have $a_2 = a_3$.  Let $P'
= P_+(T_1) = P_-(T_2)$ be the twice punctured disk cutting $E(T)$ into
$E(T_1)$ and $E(T_2)$.  Then on $\bdd E(T_1)$ there are $a_2 + a_3$
arcs of $F \cap P'$ with both endpoints on the circle $P' \cap P(T)$
while there are $a_5 + 2a_6$ such arcs on $\bdd E(T_2)$.  Hence $2a_2
= a_2 + a_3 = a_5 + 2a_6 = n + 2s$, which shows that $n = 2r$ is even,
and $a_2 = a_3 = r+s$.  (2) follows from this and the equation $a_1 +
a_2 = n = 2r$, as shown in Figure 5.2(a).

Conversely, given a set of non-negative numbers $a_k$ satisfying the
above equations, let $E_k$ be $a_k$ copies of disks of type (k).  Then
one can check that the arcs $(E_1 \cup E_2 \cup E_3) \cap P'$ are
isotopic to the arcs $(E_4 \cup E_5 \cup E_6) \cap P'$ on $P'$, hence
one can glue these together to form a surface $F$ in $E(T)$.
\qed

The curves in Figures 5.2(a) and 5.2(b) can be represented by the
weighted train tracks shown in Figure 5.2(c) and 5.2(d), respectively,
where an arc of the train track with weight $x$ represents $x$
parallel copies of that arc.  The weights on the train tracks are
calculated using the above lemma.

Recall that $T(1/3, -1/2) = T_1 + T_2$ is formed by gluing $P_+(T_1)$
to $P_-(T_2)$ using a reflection map along the $y$ axis on $R^2$.
Since the train track on these two half planes match each other under
this reflection, after gluing we obtain a special surface $F$ in
$E(T)$ with $F \cap P$ represented by the train track $\gamma$ in
Figure 5.3(a).

To simplify the diagram, we perform a counterclockwise full twist on
both the top two tangle endpoints and the bottom two tangle endpoints.
Note that this is equivalent to twisting the two lower endpoints by
four left hand half twists, so by definition it deforms the tangle
$T(1/3, -1/2)$ to the tangle $T(1/3, -1/2;\, 4)$ defined in Section 2.
The train track, after splitting along two edges, becomes that in
Figure 5.3(b).

\bigskip
\leavevmode

\centerline{\epsfbox{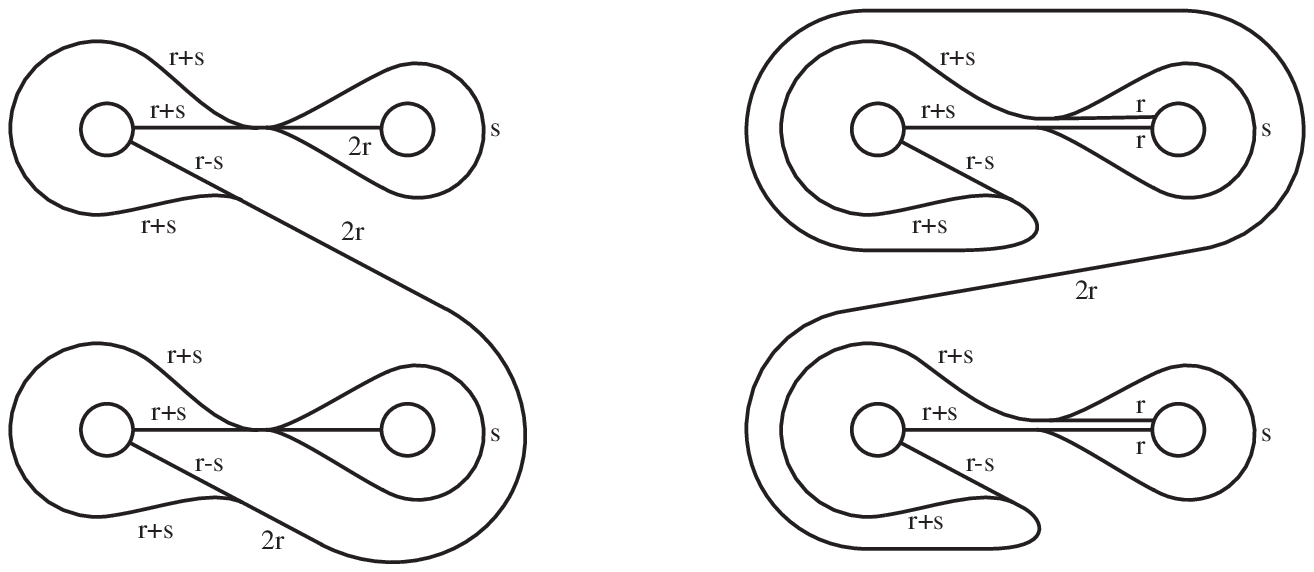}}
\bigskip
\centerline{Figure 5.3}
\bigskip

After a further splitting and isotopy, we obtain the train track
$\gamma$ in Figure 5.4(a).  Note that up to isotopy we can move the
end points of $\gamma$ around $\bdd P$, but so far we have not done
that.  By moving some endpoints of train track around $\bdd P$, we
obtain the one in Figure 5.4(b).

\bigskip
\leavevmode

\centerline{\epsfbox{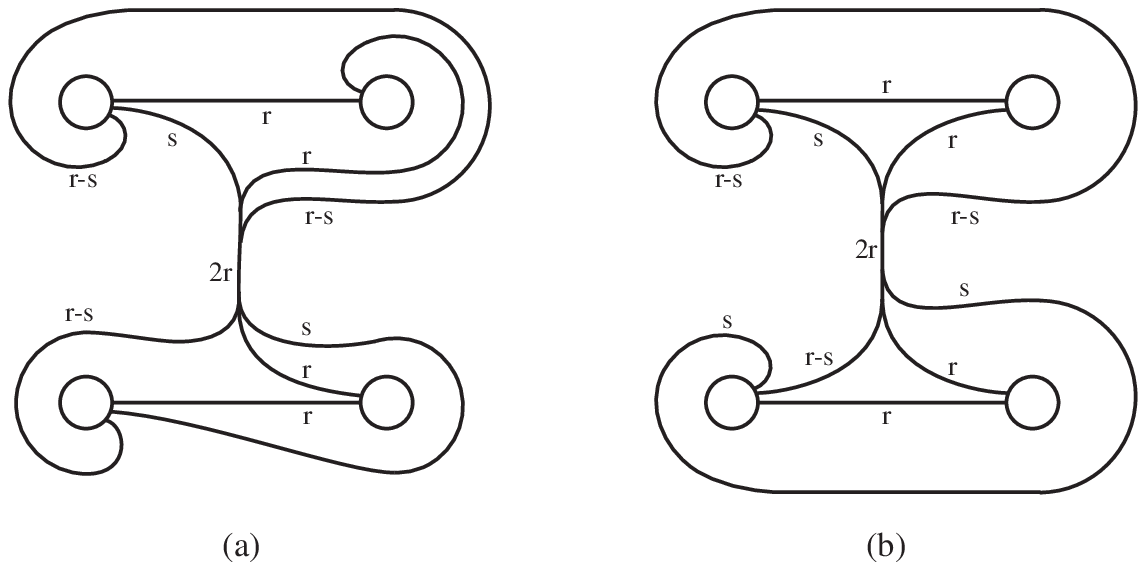}}
\bigskip
\centerline{Figure 5.4}
\bigskip

The two strings of $T = T(\pm 1/q)$ is said to be {\it consistently
  oriented\/} if they both run from the upper endpoints to the lower
endpoints or both from the lower endpoints to the upper endpoints.
For $T(1/3)$ we will always assume that its two strings are
consistently oriented.  For $T(-1/2)$, we introduce a new variable
$\epsilon$ and set $\epsilon = 1$ if the two strings of $T(-1/2)$ are
oriented consistently, and $\epsilon = -1$ otherwise.  Recall that a
surface $F$ in $E(T)$ is regular if $F\cap U(T)$ is a set of regular
curves on $U(T)$.

\begin{lemma} Let $A_k$ be a regular special disk such that $A_i \cap
  P(T)$ is the curve in Figure 2.1(k).

  (1) $\theta(A_1) = 6$, $\theta(A_2 \cup A_3)= 4$, $\theta(A_4 \cup
  A_5) = -2 \epsilon$, and $\theta(A_6) = 0$.

  (2) A special surface $F$ in $T(1/3, -1/2)$ is isotopic to a regular
  surface $F'$ such that $\bdd F' \cap P$ is represented by the
  weighted train track in Figure 5.3(a), and $\theta(F') = (5 -
  2\epsilon) n - 2s$, where $s = a_6$ is the number of type (6)
  special disks in $F$.

  (3) A special surface $F$ in $T(1/3, -1/2;\, 4)$ is isotopic to a
  regular surface $F'$ such that $F' \cap P$ is carried by the
  weighted train track in Figure 5.4(b), and $\theta(F) = (6 -
  4\epsilon) n - 2s$.
\end{lemma}

\proof (1) follows by drawing the curves $\gamma$ of $\bdd A_i$ on
$U(T)$ and counting the signed crossings where $\gamma$ passes below
$\tau$.  We omit the details.

(2) By definition $F$ is the union of $a_k$ copies of $A_k$, which can
be put in regular position.  We leave it to the reader to check that
the regular surfaces in $E(T(1/3))$ and $E(T(-1/2))$ can be combined
together to create the surface $F$ in $E(T(1/3, -1/2))$ without
creating new crossings between $F \cap U(T)$ and $\tau$.  Since $a_2 =
a_3$ and $a_4 = a_5$, we have $\theta(F) = \theta(A_1) + a_2 \,
\theta(A_2 \cup A_3) + a_4 \, \theta(A_4 \cup A_5) + a_6 \,
\theta(A_6)$.  It follows from (1) and Lemma 5.3 that
\begin{eqnarray*}
  \theta(F) & = &  6\, a_1 + 4\, a_2 - 2\, \epsilon\, a_4 \\
  & = & 6(r-s) + 4(r+s) - 2\, \epsilon\, n = (5-2\epsilon) n - 2s. 
\end{eqnarray*}

(3) By definition $T = T(1/3, -1/2;\, 4)$ is obtained from $T(1/3,
-1/2)$ by two counterclockwise full twists of the two lower endpoints
of the tangle, which deforms the surface $F''$ in (2) to a new surface
$F'$ in $E(T)$ with boundary curve represented by the train track
shown in Figure 5.4(a).  Note that after the twists each arc component
of $\alpha' = F' \cap \bdd N(\tau)$ passes below $\tau$ four more
times than $\alpha'' = F'' \cap \bdd N(\tau)$ does, two of which in
the positive direction, and the other two in positive direction if and
only if $\epsilon=-1$.  (Note that the two strings twisted have the
same orientation if and only if the two strings in the tangle
$T(-1/2)$ have opposite orientation.)  Hence $\theta(F') = \theta(F'')
+ (2-2\epsilon)n$, so by (2) we have $\theta(F') = (7 - 4\epsilon) n -
2s$.

We now perform an isotopy of $F'$ to obtain the surface $F$ whose
boundary curve on $P$ is carried by the train track in Figure 5.4(b).
To do this one needs to turn $2r = n$ endpoints of $F \cap P$ on $\bdd
P$ clockwise for an angle of almost $2\pi$.  The isotopy on each
endpoints creates one more crossing at which $\alpha'$ passes below
$\tau$, and it is in the negative direction.  Therefore we have
$\theta(F) = \theta(F') - n = (6 - 4\epsilon) n - 2s$.  \qed

The following lemma can be used to calculate the boundary slope of
a surface in $E(K)$.

\begin{lemma} Suppose $F_i$ is a regular surface in $E(T_i)$.  Let
  $\eta: \bdd B_1 = \hat R^2 \to \hat R^2 = \bdd B_2$ be the
  reflection along the $y$ axis, such that $\eta(F_1\cap P(T_1)) =
  F_2\cap P(T_2)$.  Let $(S^3, K) = (B_1, \tau_1) \cup _{\eta} (B_2,
  \tau_2)$ with orientation of $\tau_i$ induced by that of $K$, and
  let $F = F_1 \cup _{\eta} F_2$.  Suppose $F$ has $m$ boundary
  components with slope $p/q$, where $p, q$ are coprime and $q>0$.
  Then $mp = \theta(F_1) + \theta(F_2)$.  In particular, $q = 1$
  if and only if $\theta(F_1) + \theta(F_2) \equiv 0$ mod $n$, where
  $n = mq$ is the number of times $\bdd F$ intersects each meridian of
  $K$.
\end{lemma}

\proof We can shift $T_1=(B^3, \tau_1)$ to the left and $T_2=(B^3,
\tau_2)$ to the right so that $\tau_1, \tau_2$ are separated by the
$yz$-half-plane in $B^3 = R^3_-$.  Now $K\subset R^3$ is isotopic to
the knot $K'$ obtained by adding four arcs on $R^2$, each connecting
an endpoint $p_i$ of $\tau_1$ to $\eta(p_i)$ on $\tau_2$, with two
below the line $y=-1$ and the other two above the line $y=1$.  We may
also assume that near $\bdd \tau_i$ these arcs match the arcs $\alpha$
in Figure 5.1, so they are disjoint from the projection of $F_i \cap
U(T_i)$ because $F_i$ are regular.  The isotopy from $K$ to $K'$
extends to an isotopy which deforms $F$ in $E(K)$ to the surface $F'$
in $E(K')$ obtained from $F_1 \cup F_2$ by connecting their boundary
on $R^2$ by bands in the upper half space.  More explicitly, let $C =
F_1\cap R^2$ and embed $C \times I$ in $R^3_+ \cap E(K')$ so that
$C\times \{-1\} = F_1 \cap R^2$, $C\times \{1\} = F_2 \cap R^2$, and
$\bdd C \times I$ lies on $\bdd N(K')$.  Then $F' = F_1 \cup (C\times
I) \cup F_2$.  Note that $(\bdd C)\times I$ does not pass below $K'$,
so $lk(K', \bdd F')$ is the sum of the signs of crossings where $\bdd
F_i$ passes below $\tau_i$, hence $mp = lk(K, \bdd F) = lk(K', \bdd
F') = \theta(F_1) + \theta(F_2)$.  

If $q=1$ then $mp\equiv 0$ mod $n = mq$.  Conversely, if $mp \equiv 0$
mod $n = mq$ then $p$ is a multiple of $q$.  Since $p,q$ are coprime
and $q>0$, this is possibly only if $q=1$.  
\qed

We now assume that $K$ is a type II knot and $F$ is an essential
punctured torus in $E(K)$ with integer boundary slope.  By Proposition
4.4, $K$ is the union of two tangles $T_1, T_2$, each weakly
equivalent to $T(1/3, -1/2)$ or $T(-1/3, 1/2)$.  Hence up to weak
isotopy we may assume that $T_i = T(1/3, -1/2;\, 4)$ or its mirror
image $T(-1/3, 1/2; -4)$.  By Lemma 4.3 we may assume that $F_i = F
\cap E(T_i)$ is a special surface, so by Lemma 5.4(3) we may assume
that $F_i$ is regular, and $F_i \cap P(T_i)$ is represented by the
train track $\gamma_i$ in Figure 5.4(b) if $T_i = T(1/3, -1/2;\, 4)$,
and its reflection along the $y$ axis if $T_i = T(-1/3, 1/2; -4)$.
Note that the gluing map $\eta: \bdd B_2 \to \bdd B_1$ could be any
orientation reversing map which maps $P(T_2)$ to $P(T_1)$ and $F_1
\cap P(T_1)$ to $F_2\cap P(T_2)$.

As in Lemma 5.3, let $n$ be the number of times $F$ intersects a
meridian of $K$ on $N(K)$, $r = n/2$, and let $s_i$ be the number of
type (6) special disks in $F_i$.  There are five possible ways to
split $\gamma_i$, according to the values of $s_i$.  

(1) $r-s_i > s_i > 0$; 

(2) $s_i > r-s_i > 0$; 

(3) $s_i = 0$; 

(4) $r-s_i = 0$; and 

(5) $r-s_i = s_i$.  

One can check that for $T_i = T(1/3, -1/2;\, 4)$ the train track in
Figure 5.4(b) splits to $\gamma_i$ in Figure 5.5(1)--(5),
respectively, where $s = s_i$.  When $T_i = T(-1/3,1/2; -4)$,
$\gamma_i$ is the reflection of those in the figure along the
$y$-axis.  We say that $\gamma_i$ is of type (k) if it is the one in
Figure 5.5(k).

\bigskip
\leavevmode

\centerline{\epsfbox{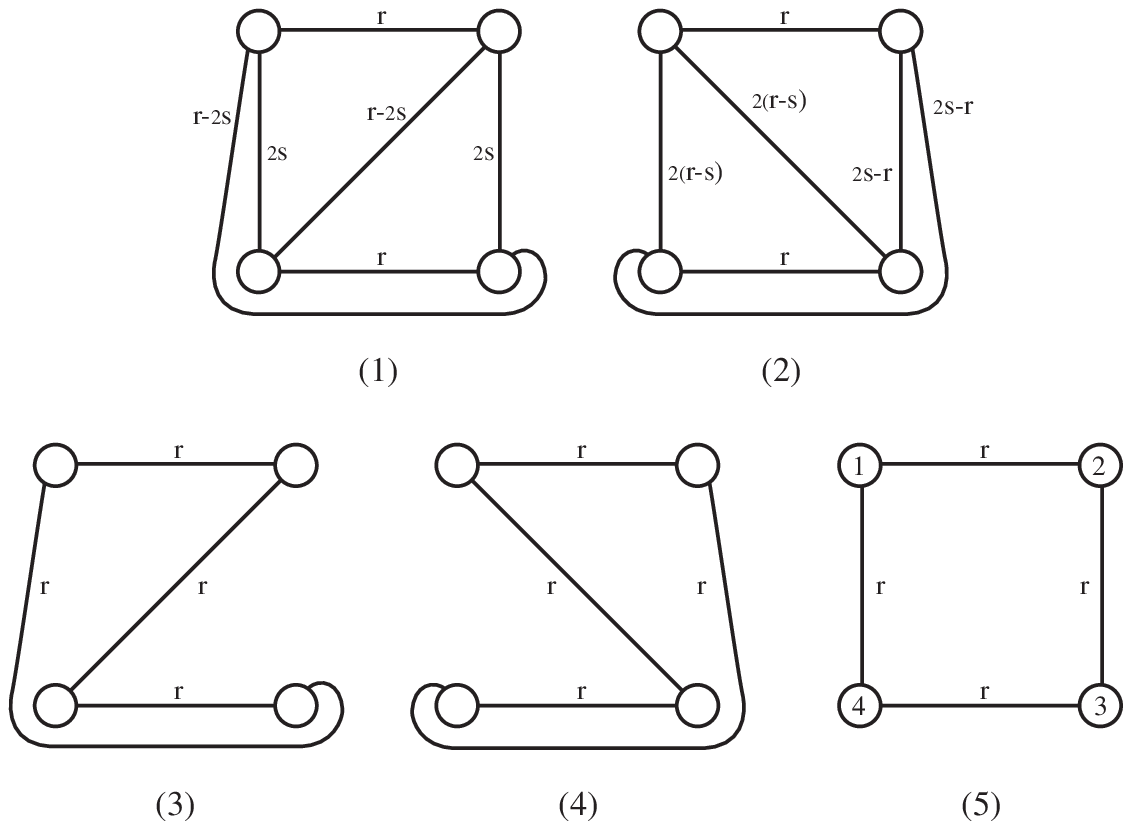}}
\bigskip
\centerline{Figure 5.5}
\bigskip

\begin{lemma} 
For each $i$, $s_i \in \{0,\, r,\, r/2\}$.  Hence $\gamma_i$ is not of
type (1) or (2).
\end{lemma}

\proof We need to show that $\gamma_i$ cannot be of type (1) or (2).
Because of symmetry we may assume without loss of generality that $T_1
= T(1/3, -1/2, 2)$, and $\gamma_1$ is of type (1).  By Lemma 5.4(3) we
have $\theta(F_1) = (6 - 4\epsilon) n - 2s_1$.

Since the graph in Figure 5.5(1) is not homeomorphic to those in
Figure 5.5(3)--(5), we see that $\gamma_2$ must be of type (1) or (2).
Let $\eta: \bdd B_2 \to \bdd B_1$ be the gluing map, which is
orientation reversing.  Since the two horizontal edges of $\gamma_i$
have higher weights than the other edges, $\eta$ must map horizontal
edges to horizontal edges.  Without loss of generality we may assume
that $\eta$ maps the upper horizontal line of $\gamma_2$ to the upper
horizontal line of $\gamma_1$ by a reflection along the vertical axis.
Note that this completely determine $\eta$ on $\gamma_2$.

First assume that $T_2 = T(1/3, -1/2;\, 4)$.  If $\gamma_2$ is the one
in Figure 5.5(2) then $\eta$ is simply a reflection along the
vertical line, hence it maps the two right vertices of $\gamma_2$ to
the two left vertices of $\gamma_1$, but since the two left (right)
vertices of $\gamma_i$ belong to the same component of $\tau_i$, $K$
would be a link of two components, which is a contradiction.  Therefore
$\gamma_2$ must be the graph in Figure 5.5(1), which is redrawn in
Figure 5.6(1).  One can modify the graph by turning the lower
horizontal edge clockwise by a half twist to obtain the one in Figure
5.6(2), then isotope some of the edge endpoints at the two lower
boundary components of $P$ around to obtain the graph in Figure
5.6(3).  The map $\eta$ is the composition of this isotopy followed
by a reflection along a vertical line.

\bigskip
\leavevmode

\centerline{\epsfbox{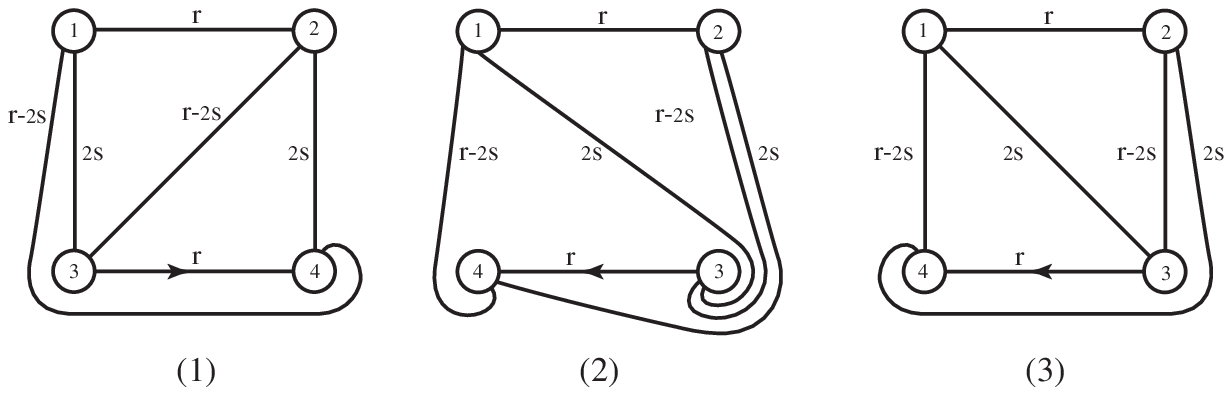}}
\bigskip
\centerline{Figure 5.6}
\bigskip

Note that the isotopies above are not regular isotopies.  They have
changed the relative framing of the surface $F_2$.  By Lemma 5.4(3)
the framing of $F_2$ with boundary graph in Figure 5.6(1) is given by
$(6 - 4\epsilon) n - 2s_2$.  After twisting the two lower vertices of
the graph by a half twist, each boundary arc of $F_2$ on the tubes
$\bdd E(T_2) - P$ passes below the part of $\tau$ near the vertex $3$
in the figure once in the negative direction, but does not pass below
the other string of $\tau$, hence the new framing is $(4 - 4\epsilon)
n - 2s_2$.  The isotopy from Figure 5.6(2) to Figure 5.6(3) moves $r$
edge endpoints clockwise and another $r$ edge endpoints
counterclockwise around vertex $3$ and $4$ respectively, hence will
not change the relative framing.  Therefore we have $\theta(F_2) = (4
- 4\epsilon) n - 2s_2$ for the surface $F_2$ corresponding to Figure
5.6(3).

In order to glue $F_2$ to $F_1$ by $\eta$, the weight of the left
vertical edge in Figure 5.6(3) must match the weight of the right
vertical edge of Figure 5.5(1), hence we have $r - 2s_2 = 2s_1$.  Thus
$\theta(F_1) = (6-4\epsilon)n - (r - 2s_2)$, so 
\begin{eqnarray*}
\theta(F_1) + \theta(F_2) & = & [(6-4\epsilon)n - (r-2s_2)] +
[(4-4\epsilon)n - 2s_2] \\
& = & (10-8\epsilon)n - r \equiv r \mbox{\qquad mod $n$}.
\end{eqnarray*}

Since $r = n/2$, this is a contradiction to Lemma 5.5 and the
assumption that $F$ has integer boundary slope.  \qed

\begin{lemma} 
$s_1 = s_2 = r/2$, so both $\gamma_i$ are of type (5).
\end{lemma}

\proof   Note that if $\gamma_i$ is of type (3) or (4) then the
endpoints of a string of $\tau_i$ is separated by those of the other
string on $\gamma_i$, which is a circle.  Hence if both $\gamma_i$ are
of type (3) or (4) then $K = \tau_1 \cup \tau_2$ would be a link of
two components, which is a contradiction.  Therefore at least one
$\gamma_i$, say $\gamma_1$, is of type (5).  If the result is not true
then $\gamma_2$ must be of type (3) or (4).  We assume it is of type
(3).  The other case is similar.

By the same proof as in that of Lemma 5.6, we can isotope $\gamma_2$
by twisting the lower level edge of $\gamma_2$ by a half twist,
followed by an isotopy which moves some endpoints of $\gamma_2$ around
$\bdd P$, to change $\gamma_2$ to a graph of type (5).  As in the
proof of Lemma 5.6, this will not change $\theta(F_2)$ mod $n$.  By
Lemma 5.4(3) we have $\theta(F_1) + \theta(F_2) \equiv -2s_1 -2s_2$ mod
$n$.  Since $\gamma_1$ is of type (5), we have $s_1 = 0$, and since
$\gamma_2$ is of type (3), we have $2s_2 = r = n/2$.  Hence
$\theta(F_1) + \theta(F_2) \equiv r$ mod $n$, which by Lemma 5.5
implies that the boundary slope of the punctured torus $F$ is not an
integer slope, a contradiction.  \qed

\begin{prop} Let $K$ be a type II knot, and let $F$ be an essential
  punctured torus in $E(K)$ with integer boundary slope $\delta$.
  Then $(K, \delta) = (K_1, 3)$, $(K_2, 0)$ or $(K_3, -3)$, where
  $K_i$ are the knots defined in Definition 2.3.
\end{prop}

\proof By Proposition 4.4 we have $(S^3, K) = T_1 \cup T_2$, where
each $T_1$ is weakly equivalent to $T(1/3, 1/2)$ or $T(-1/3, 1/2)$.
Up to weak equivalence we may assume that each $T_i$ is either $T(1/3,
-1/2;\, 4)$ or $T(-1/3, 1/2; -4)$.

Denote by $u_1, ..., u_4$ the four disks $\bdd B_i - Int P(T_i)$,
which we will consider as fat vertices.  By Lemma 5.7 the train track
$\gamma_1$ is of type (5), so it is a cycle containing those four
vertices, labeled in cyclic order, as shown in Figure 5.5(5).
Similarly, $\gamma_2$ contains the vertices $v_1, ..., v_4$ in the
same order.  Orient $\gamma_i$ clockwise.  Let $\eta: \bdd B_1 \to
\bdd B_2$ be the gluing map.  Then up to isotopy $\eta$ is determined
by its restriction on $\gamma_1$, which in tern is determined by the
image of $u_1$ and whether $\eta|_{\gamma_1}$ is orientation
preserving or not.  Note that although $\eta$ must be orientation
reversing on $\bdd B_1$, it may map the disk inside of $\gamma_1$ to
the disk outside of $\gamma_2$, so $\eta|_{\gamma_1}$ could be
orientation preserving.

Since $K$ is a knot, the endpoints of a string of $\tau_2$ must be
mapped to endpoints of different strings of $\tau_2$, which excludes
four possible $\eta$.  Also, if $\tau_i$ is considered as lying in a
pillow case $B_i$, then from Figure 2.2(a) one can see that a $\pi$
rotation along a horizontal axis will preserve the tangle.  One can
now easily check that all the four possible choices of $\eta$ give
rise to the same knot, so we may assume that $\eta$ is obtained by
rotation $\gamma_2$ counterclockwise by an angle of $\pi/2$ followed
by a reflection along a vertical line, as described in Definition 2.3.
Therefore $K$ is one of the three knots in the statement.

The surface $F$ is cut into $F_i$ in $E(T_i)$.  By Lemma 5.2(3), the
$\pi/2$ rotation above will not change $\theta(F_1)$, so by Lemma 5.5
the boundary slope of $F$ is given by $(\theta(F_1) + \theta(F_2))/n$.
Examining the orientation of the strings in $\tau_i$ we see that they
are consistently oriented in the tangle $T(\pm 1/2)$, hence $\epsilon
= 1$ for both $T_i$.  Since $\gamma_i$ are of type (5), by definition
we have $2s_i = r$, so by Lemma 5.4(3) we have $\theta(F_i) = (6 - 4
\epsilon)n - 2s_i = 2n - r$ if $T_i = T(1/3, -1/2;\, 4)$, and
$\theta(F_i) = -2n + r$ if $T_i = T(-1/3, 1/2;\, -4)$.  Hence if $K =
K_1 = T(1/3, -1/2;\, 4) \cup_{\eta} T(1/3, -1/2;\, 4)$ then by Lemma
5.5 the boundary slope of $T$ is $[(2n-r)+(2n-r)]/n = 3$.  Similarly
for the other two cases.  \qed

\section{Toroidal surgery}

Let $(K, \delta)$ be one of the three pairs described in Theorem 1.1.
In this section we will show that there is a punctured torus $F$ in
$E(K)$ with boundary slope $\delta$, and the torus $\hat F$ obtained
by capping off the boundary components of $F$ with meridional disks in
the Dehn filling solid is indeed an essential torus in the surgered
manifold $K(\delta)$.

Let $T = T(1/3, 1/2) = T_1 \cup T_2$, where $T_1 = T(1/3)$ and $T_2 =
T(1/2)$.  

\begin{lemma}
  A special surface $Q$ in $E(T)$ is incompressible and
  $P$-incompressible.
\end{lemma}

\proof By considering a component if necessary we may assume that $Q$
is connected.  By Lemma 3.3 each special disk has zero index, so by
the Additivity Lemma 3.2 $Q$ also has zero index and hence is either
an annulus disjoint from $U_+(T)$ or a disk intersecting $U_+(T)$ in
two arcs.  If $Q$ is a disk then it is automatically incompressible,
and a $P$-compression will produce two disks $D_i$, each intersecting
$U_+(T)$ at a single arc, which contradicts Lemma 3.4.

Now assume $Q$ is an annulus.  Let $P' = P_+(T_1) = P_-(T_2)$ be the
twice punctured disk cutting $E(T)$ into $E(T_1)$ and $E(T_2)$.  By
definition $P'$ cuts $Q$ into a set of special disks, each of which is
an essential bigon in the sense that it intersects $P'$ in two arcs,
and an arc on the bigon with one endpoint on each of these arcs is not
rel $\bdd$ homotopic to an arc on $P'$.  Using this and the fact that
$P'$ is incompressible one can easily show, by an innermost circle
outermost arc argument, that $Q$ is incompressible in $E(T)$.  To show
it is $P$-incompressible, one need only show that there is no
$\bdd$-compressing disk $D$ of $Q$ in $E(T)$ disjoint from $U_+(T)$.
Suppose $\bdd D = \alpha \cup \beta$, where $\alpha$ is an essential
arc on $Q$ and $\beta \subset \bdd E(T) - U_+(T)$.  Since $P'$ is
incompressible, we may assume $D \cap P'$ has no circle components;
since $P' \cap Q$ is a set of essential arcs on $Q$, $\alpha$ is
isotopic to a component of $P' \cap Q$, so by an isotopy of $D$ we may
also assume that $\alpha \cap P' = \emptyset$.  Hence $D \cap P'$ is a
set of arcs with endpoints on $\beta$.  Choose $D$ so that $|D \cap
P'|$ is minimal.  If $D \cap P' = \emptyset$ then one can use $D$ to
$\bdd$-compress a special disk to produce a pair of disks $D_j$ in
some $E(T_i)$ with $\bdd D_j$ the union of an essential arc on $P'$
and another arc on $\bdd E(T_i) - P'\cup U_+(T)$, which will lead to a
contradiction to Lemma 3.3.  Now an outermost arc $\gamma$ cuts off a
disk $D'$ from $D$ which lies in one of the $E(T_i)$.  Using Lemma 3.3
one can show that $\bdd D'$ is trivial on $\bdd E(T_i)$, so it cuts
off a 3-ball which can be used to reduce $|D\cap P'|$, contradicting
its minimality.  \qed

By definition the knot $K$ is the union of two tangles $T_1, T_2$,
each of which is either $T(1/3, -1/2;\, 4)$ or its mirror image
$T(-1/3, 1/2;\, -4)$.  Let $n=4$, $s=1$, and define $a_k$ by $a_1 =
a_6 = 1$, $a_2=a_3 = 3$, and $a_4=a_5=4$.  One can check that these
numbers satisfy the equations in Lemma 5.3.  Since $T(1/3, -1/2; 4)$
is weakly equivalent to $T(1/3, 1/2)$, by Lemma 5.3 there is a special
surface $F_i$ in $E(T_i)$ which is the union of $a_k$ copies of
special disks of type (k) for $k=1,...,6$.  By Lemma 5.4(3) the curves
$F_i \cap P(T)$ is represented by the train track in Figure 5.4(b),
which splits to the one in Figure 5.5(5) because $r=2=2s$.  Similarly
for $T_i = T(-1/3, 1/2; -4)$.  The graph in Figure 5.5(5) is
preserved by the gluing map $\eta: \bdd B_1 \to \bdd B_2$, which by
definition is a $\pi/2$ rotation followed by a reflection along the
vertical line.  It follows that $F_1 \cup_{\eta} F_2$ form a surface
$F$ in $E(K)$.

\begin{lemma} $F$ is an essential punctured torus in $K$ with $\bdd F$
  consisting of $4$ circles of slope $\delta$.
\end{lemma}

\proof The boundary slope $\delta$ of $F$ is calculated in the proof
of Proposition 5.8.  In particular, $\delta$ is an integer in all
three cases of $K$.  Since $n =4$ in the construction, we see that
$|\bdd F| = 4$.

By Lemma 6.1 each $F_i$ is incompressible and $P$-incompressible, and
by Lemma 3.4 $P = P(T_1) = P(T_2)$ is also incompressible.  Thus by an
innermost circle outermost argument one can show that $F = F_1 \cup
F_2$ is incompressible in $E(K)$.  Since $F$ has four boundary
components, this implies that $F$ is also $\bdd$-incompressible as
otherwise two copies of a boundary compression disk and the annulus on
$\bdd N(K)$ bounded by two components of $\bdd F$ would form a
compressing disk of $F$.  \qed

\begin{lemma} Let $M$ be a handlebody of genus $3$, let $c_1, c_2$ be
  a pair of curves on $\bdd M$, and let $M'$ be the manifold obtained
  by attaching two $2$-handles along $c_1$ and $c_2$.  If $\bdd M -
  c_i$ is compressible for $i=1,2$, and $\bdd M - c_1 \cup c_2$ is
  incompressible, then $\bdd M'$ is incompressible.
\end{lemma}

\proof This is a standard application of the Handle Addition Lemma
[Ja].  Denote by $M_1$ the manifold obtained from $M$ by attaching a
2-handle along $c_1$.  By assumption there is a compressing disk $D$
of $\bdd M$ which is disjoint from $c_1$.  If $D$ is separating then
one component $H$ of $M|D$ is a handlebody disjoint from $c_1$, so we
may re-choose $D$ to be a non-separating disk in $H$.  Thus after
attaching a 2-handle to $M$ along $c_1$, $D$ is still a compressing
disk of $\bdd M_1$, so $\bdd M_1$ is compressible.  On the other hand,
since $\bdd M - c_2$ is compressible while $(\bdd M - c_2) - c_1$ is
incompressible, by the Handle Addition Lemma applied to the pair
$(M-c_2, c_1)$ we see that the surface $\bdd M_1 - c_2$ is
incompressible in $M_1$.  Now since $\bdd M_1$ is compressible while
$\bdd M_1 - c_2$ is incompressible, we may apply the Handle Addition
Lemma again to conclude that after attaching a 2-handle to $M_1$ along
$c_2$, the boundary of the resulting manifold $M'$ is incompressible.
\qed

Let $F_0$ be a separating surface in a 3-manifold $M$ with $\bdd F_0 =
c_1 \cup ... \cup c_n \cup c'_n \cup ... \cup c'_1$, lying
successively on a torus component $R$ of $\bdd M$.  Let $A_1$ be the
component of $R|C$ bounded by $c_1 \cup c'_1$, and let $A_k$ be the
annulus on $R$ which is bounded by $c_k \cup c'_k$ and contains $A_1$.
Starting with $F_0$, one one can construct a sequence of surfaces
$F'_k$ and $F_k$ by adding the annulus $A_k$ to $F_{k-1}$ to obtain
$F'_k$ and then pushing the $A_k$ part of $F'_k$ off $R$ to obtain
$F_k$.  The surfaces $F_k$ are said to be obtained from $F$ by {\it
  successively tubing through $A_1$}.  The following lemma is probably
due to Gordon.

\begin{lemma} Suppose $F_0$ is a connected separating incompressible
  surface in a 3-manifold $M$ with $\bdd F_0$ on a torus component $R$
  of $\bdd M$.  Let $M'_0, M''_0$ be the components of $M|F_0$, and
  let $A_1$ be an annulus component of $\bdd M'_0 - \Int F_0$.  If
  $F'_1 = F_0 \cup A_1$ is incompressible in $M'_0$, then the surfaces
  $F_k$ obtained by successively tubing $F$ through $A$ are all
  incompressible in $M$.
\end{lemma}

\proof We use the above notation, and let $M'_k, M''_k$ be the
components of $M|F_k$.  Note that $F_k$ are all connected and
separating, and $A_{k+1}$ is a component of $\bdd M'_k - \Int F_k$ or
$\bdd M''_k - \Int F_k$.  Hence by induction we need only show that
(i) $F_1$ is incompressible, and (ii) if $n>1$ then $F'_2$ is
also incompressible in $M''_1$.

Since $F_1$ is obtained by pushing the $A_1$ part of $F'_1 = F_0 \cup
A_1$ into the interior of $M$, the components $M'_1$ of $M|F_1$
containing $A_1$ is homeomorphic to $M'_0$, with $F_1 \subset \bdd
M'_1$ identified to $F_0 \cup A_1 \subset \bdd M'_0$, so by assumption
$F_1$ is incompressible in $M'_1$.  When $n=1$ the other component
$M''_1$ is obtained by attaching a collar $R \times I$ to $M''_0$
along the annulus $A'_1 = R - \Int A_1$.  It is clear that $A'_1$ is
incompressible in $M''_1$.  If it is also $\bdd$-incompressible then
an innermost circle outermost arc argument shows that $F_1$ is
incompressible in $M''_1$, and if it is $\bdd$-compressible then $F_0$
must be an annulus which is parallel to $A'_1$ after cutting off some
possible summands, so $M''_1$ is essentially a product $R\times I$ and
hence $F_1$ is also incompressible.  When $n>1$ $M''_1$ is obtained by
attaching $A_2 \times I$ to $M''_0$ along the two annuli $P$ and $P'$
bounded by $c_1 \cup c_2$ and $c'_1 \cup c'_2$, respectively.  The
incompressibility of $F_0$ implies that these annuli are
incompressible.  Also, there is no disk $D$ in $M''_1$ intersecting $P
\cup P'$ at a single essential arc as otherwise the frontier of $N(D
\cup P \cup P')$ would contain a compressing disk of $F_0$.  One can
now apply an innermost circle outermost arc argument to show that both
$F_1$ and $F'_2$ are incompressible in $M''_1$.  \qed

\medskip
\noindent {\bf Proof of Theorem 1.1.\/} If $K$ is not a type II knot
or if $\delta$ is not an integer slope then by [Wu2, Theorems 3.6 and
4.4] $K(\delta)$ is Haken and hyperbolic, so we assume that $K$ is a
type II knot and $\delta$ is an integer slope.  Write $(S^3, K) = T_1
\cup T_2$, with $T_i = T(p_i/q_i, 1/2)$.  By [Wu2, Theorem 2.3] $E(K)$
contains an essential branched surface ${\cal B}$ which remains
essential in $K(\delta)$, hence by [GO] $K(\delta)$ is irreducible.
Also, by the construction in the proof of [Wu2, Theorem 2.3] the
exterior of the ${\cal B}$ is the disjoint union of $E(T_1)$ and
$E(T_2)$, with vertical surface $U_+(T_i)$ on $\bdd E(T_i)$.

We claim that $E(T_i)$ is not an $I$-bundle with $U_+(T_i)$ as a
vertical annulus.  If this is false that after attaching a 2-handle to
$U_+(T_i)$ the resulting manifold $M_i$ would be an $I$-bundle over a
closed surface, which must be a Klein bottle because $\bdd M_i$ is a
torus.  Since $M_i$ is the exterior of a trefoil knot in $S^3$, this
would imply that there is a Klein bottle embedded in $S^3$, which is
absurd.

By [Br] if a small Seifert fiber space contains an essential
lamination then its exterior is an $I$-bundle, hence the above implies
that $K(\delta)$ cannot be a small Seifert fiber space.  Therefore
$K(\delta)$ is exceptional if and only if it is toroidal.  By
Proposition 5.8 $K(\delta)$ is toroidal only if $(K,\delta)$ is one of
the three pairs listed, so we need only show that $K(\delta)$ is
indeed a toroidal manifold for each of those pairs.

Let $F$ be the surface constructed before Lemma 6.2.  By Lemma 6.2 $F$ is an
essential punctured torus with $\bdd F$ consisting of four circles of
slope $\delta$ given in Theorem 1.1.

Let $F_i = F\cap E(T_i)$, $P = P(T_i)$, and $P_i = E(T_{iq}) \cap
E(T_{i2})$.  Note that the special disks $F \cap E(T_{ij})$ in
$E(T_{ij})$ cuts $E(T_{ij})$ into a set of 3-balls $B_k$.  The arcs
$P_i \cap F_i$ cut $P_i$ into a set of disks $D_r$, so the manifold
$E(T_i)|F_i$ is obtained by gluing the $B_k$'s along the $D_r$'s, and
hence is a set of handlebodies $H_k$.  Also, $F\cap P$ cuts $P$ into a
set of disks $D'_r$, so $E(K)$ is obtained from the $H_k$'s by gluing
along the $D'_r$'s, and hence is also a set of handlebodies.  Since
$F$ is a punctured torus with four boundary components, it is
separating in $E(K)$, so $E(K)|F$ has two components $M_1, M_2$, and
$\bdd M_i$ is the union of $F$ with two annuli on $\bdd N(K)$ and
hence is a surface of genus $3$.  It follows that each $M_i$ is a
handlebody of genus $3$.

Let $A_1, A_2$ be the two annuli $\bdd M_i -\Int (F)$, and let $c_j$
be the core of $A_j$.  Then by Lemma 6.2 the surface $\bdd M_i - c_1
\cup c_2$, which is homotopic to $F$ on $\bdd M_i$, is incompressible.
Note that $\bdd M_i - c_1$ is homotopic to the surface $F \cup A_2$
obtained from $F$ by tubing along $A_2$.  If this is incompressible
then by Lemma 6.4 the closed surface $F'$ obtained from $F$ by
successively tubing through $A_2$ is incompressible.  From the
construction one can see that $F'$ has coannular slope $\delta$ on
$\bdd N(K)$ in the sense that there is an incompressible annulus $A'$
which has interior disjoint from $F'$, with one boundary component on
$F'$ and the other on $\bdd N(K)$ with slope $\delta$.  It is easy to
see that any embedded essential surface in an irreducible 3-manifold
has at most one coannular slope on a torus boundary component, hence
there is no disk in $S^3$ with boundary on $F'$ which intersects $K$
exactly once, so $F'$ is $K$-essential in the sense of [Wu2].  By
[Wu2, Lemma 4.7] there is no such closed surface in the exterior of a
type II knot, which is a contradiction.

Let $\hat M_1, \hat M_2$ be the two components of $K(\delta) | \hat
F$.  Then $\hat M_i$ is obtained from $M_i$ by attaching two 2-handles
along the curves $c_1, c_2 \subset \bdd M_i$.  We have shown above
that $M_i$ is a handlebody of genus $3$, $\bdd M_i - c_j$ is
compressible, and $\bdd M_i - c_1 \cup c_2$ is incompressible.
Therefore by Lemma 6.3 the surface $\hat F = \bdd \hat M_i$ is
incompressible in $\hat M_i$.  Since this is true for $i=1,2$, it
follows that $\hat F$ is an incompressible torus in $K(\delta)$.  \qed

\bigskip

\noindent
Department of Mathematics,  University of Iowa,  Iowa City, IA 52242
\\
Email: {\it wu@math.uiowa.edu}

\enddocument